\documentclass[12pt]{article}

\usepackage[cp1251]{inputenc}
\usepackage[T2A]{fontenc}
\usepackage[english,russian,ukrainian]{babel}
\usepackage{amsmath,amsfonts,amssymb,amsthm,cite}
\usepackage{geometry}
\newtheorem{theorem}{Теорема}
\newtheorem*{theorem-Personal}{Теорема Кигурадзе}

\newcommand{\nonprint}[1]{}

\begin{document}
\selectlanguage{russian}
\begin{flushleft}
\large
\textbf{В.А. Михайлец$^{1,2}$, О.Б. Пелехата$^{2}$, Н.В. Рева$^{2}$}
\medskip
\normalsize

$^{1}$Институт математики НАН Украины, Киев

{\it E-mail}: mikhailets@imath.kiev.ua

$^{2}$Национальный технический университет Украины ``КПИ имени Игоря Сикорского'', Киев

{\it E-mail}: 	o.pelehata-2017@kpi.ua, reva$\underline{ }$nadiia@ukr.net
\medskip

\Large

\textbf{О теореме Кигурадзе для линейных краевых задач}

\end{flushleft}

\normalsize

\it В работе исследуется предельное поведение решений неоднородных краевых задач для систем линейных обыкновенных дифференциальных уравнений на конечном интервале. Получено обобщение теоремы И.Т.~Кигурадзе (1987) о предельном переходе. \rm

\medskip

\it\textbf{Ключевые слова:} \rm Система обыкновенных дифференциальных уравнений, линейная краевая задача, предельный переход.

\bigskip

\nonprint{\selectlanguage{english}

\begin{flushleft}

\medskip
\large
\textbf{V.A. Mikhailets$^{1,2}$, O.B. Pelekhata$^{2}$, N.V. Reva$^{2}$}
\medskip
\normalsize

$^{1}$Institute of Mathematics of the National Academy of Sciences of Ukraine, Kyiv

{\it E-mail}: mikhailets@imath.kiev.ua

\medskip

$^{2}$National Technical University of Ukraine ``Igor Sikorsky Kyiv Polytechnic Institute'', Kyiv

{\it E-mail}: o.pelehata-2017@kpi.ua, reva$\underline{ }$nadiia@ukr.net
\medskip

\Large

\textbf{On Kiguradze theorem for linear boundary value problems}

\end{flushleft}

\normalsize

\it We investigate the limiting behavior of solutions of nonhomogeneous boundary value problems for the systems of linear ordinary differential equations. The generalization of Kiguradze theorem (1987) on passage to the limit is obtained. \rm

\medskip

\it\textbf{Keywords:} \rm system of ordinary differential equations, linear boundary value problem, passage to the limit.

\bigskip
}

\selectlanguage{russian}

Вопросы предельного перехода в системах дифференциальных уравнений встречаются во многих задачах теоретического и прикладного характера. Наиболее полно они исследованы применительно к решениям задачи Коши для систем дифференциальных уравнений первого порядка \cite{Reid,Opial,Levin-A,Levin,Nguen}. Более сложный случай линейных краевых задач изучался в работах И.Т.~Кигурадзе \cite{Kig,Kigu} и его последователей \cite{Kodliuk,Gnyp,Mikhailets}.

Рассмотрим на конечном интервале $ (a,b) $ систему $ {m \in \mathbb{N}} $ линейных дифференциальных уравнений первого порядка
 \begin{equation} \label{r1}
 y'(t)+A(t)y(t)=f(t)
 \end{equation}
 с общими неоднородными краевыми условиями
 \begin{equation} \label{r2}
By=c,
\end{equation}
где линейный непрерывный оператор
$$B\colon \ C([a,b];\mathbb{C}^m) \rightarrow {\mathbb{C}^m}.$$
Предполагается, что матрица-функция $A(\cdot)\in L([a,b];\mathbb{C}^{m\times m}),$ вектор-функция $f(\cdot)\in L([a,b];\mathbb{C}^{m}),$ а вектор $c\in \mathbb{C}^m$.

Под решением системы дифференциальных уравнений \eqref{r1} понимается абсолютно непрерывная на отрезке $[a,b]$ вектор-функция $y(\cdot)$, которая удовлетворяет равенству~\eqref{r1} почти всюду. Неоднородное краевое условие~\eqref{r2} корректно определено на решениях системы дифференциальных уравнений~\eqref{r1} и охватывает все классические виды краевых условий. Как известно (см., например, \cite{Kig}) краевая задача \eqref{r1}--\eqref{r2} является фредгольмовой. Поэтому для однозначной всюду разрешимости этой задачи необходимо и достаточно, чтобы однородная краевая задача имела только тривиальное решение.

Пусть теперь наряду с задачей \eqref{r1}--\eqref{r2} задана последовательность неоднородных краевых задач
\begin{equation}
 \label{r3}
 y'_n(t)+A_n(t)y_n(t)=f_n(t)
 \end{equation}
с краевыми условиями вида
\begin{equation}
 \label{r4}
B_ny_n=c_n,
\end{equation}
где матрицы-функции $A_n(\cdot)$, операторы $B_n$, вектор-функции $f_n(\cdot)$ и векторы $c_n$ удовлетворяют приведенным выше для задачи \eqref{r1}--\eqref{r2} условиям. Пусть решение $y(\cdot)$ задачи \eqref{r1}--\eqref{r2} и решения $y_n(\cdot)$ задач \eqref{r3}--\eqref{r4} существуют и однозначно определены. Тогда представляет интерес вопрос о том, когда при $n\rightarrow \infty$
\begin{equation}
\label{5}
\|y(\cdot)-y_n(\cdot)\|_\infty \rightarrow 0,
\end{equation}
где $\|\cdot\|_\infty$ --- sup-норма на отрезке $[a,b]$.

По-видимому, впервые этот вопрос был поставлен и исследован И.Т.~Кигурадзе~\cite{Kigu}. При этом предполагалось, что все функции в задаче являются вещественными.

Введем некоторые обозначения, необходимые для формулировок утверждений в удобной для нас форме. Положим
$$R_{A_n}(\cdot):=A_n(\cdot)-A(\cdot)\in L([a,b];\mathbb{C}^{m\times m}),$$
$$F(\cdot):=\left(
            \begin{array}{cccc}
              f_1(\cdot) & 0 & \ldots & 0 \\
              f_2(\cdot) & 0 & \ldots & 0 \\
              \vdots & \vdots & \ddots & \vdots \\
              f_m(\cdot) & 0 & \ldots & 0 \\
            \end{array}
          \right)\in L([a,b];\mathbb{C}^{m\times m}),          $$
$$F_n(\cdot):=\left(
            \begin{array}{cccc}
              f_{1n}(\cdot) & 0 & \ldots & 0 \\
              f_{2n}(\cdot) & 0 & \ldots & 0 \\
              \vdots & \vdots & \ddots & \vdots \\
              f_{mn}(\cdot) & 0 & \ldots & 0 \\
            \end{array}
          \right)\in L([a,b];\mathbb{C}^{m\times m}),$$
$$ R_{F_n}(\cdot)=F_n(\cdot)-F(\cdot), $$
$$R^{\vee}_{F_n}(t):=\int_a^tR_{F_n}(s){\rm d}s, \quad  R_{A_n}^{\vee}(t):=\int_a^tR_{A_n}(s){\rm d}s. $$

\begin{theorem-Personal}
Пусть выполнены условия:
\begin{enumerate}
\item[{\rm (0)}] Однородная краевая задача \eqref{r1}--\eqref{r2} имеет только тривиальное решение,

\item[{\rm (I)}] $ \|R^{\vee}_{A_n}\|_\infty \rightarrow 0$, $n\rightarrow \infty$,

\item[{\rm (II)}] $\|R_{A_n}(\cdot)\|_1=O(1)$, $n\rightarrow \infty$,

\item[{\rm (III)}] $B_ny\rightarrow By$, \quad $y\in C([a,b]$; $\mathbb{C}^m)$, $n\rightarrow \infty$.
\end{enumerate}
Тогда для достаточно больших $n$ задача \eqref{r3}--\eqref{r4} однозначно разрешима. Если кроме того выполнены условия на правые части задач
\begin{enumerate}
\item[{\rm (IV)}] $c_n\rightarrow c$, $n\rightarrow \infty$,

\item[{\rm (V)}] $\|R_{F_n}^{\vee}(\cdot)\|_\infty \rightarrow 0$, $n\rightarrow \infty$,
\end{enumerate}
то единственные решения задач \eqref{r3}--\eqref{r4} удовлетворяют предельному равенству~\eqref{5}.
\end{theorem-Personal}

Здесь и всюду дальше $\|\cdot\|_1$ --- норма в пространстве Лебега $L_1$ на отрезке $[a,b]$. Примеры показывают, что в теореме Кигурадзе все условия существенные и не одно из них нельзя отбросить. Однако, как выяснилось, некоторые из них можно значительно ослабить.

Обозначим через $\mathcal{M}^{m}:=\mathcal{M}(a,b;m)$, $m\in \mathbb{N}$ класс последовательностей матриц-функций $R_n(\cdot)\colon\mathbb{N}\rightarrow L([a,b];\mathbb{C}^{m\times m})$, для которых решение $Z_n(\cdot)$ задачи Коши
$$Z'_n(\cdot)+R_n(\cdot)Z_n(\cdot)=O, \quad Z_n(a)=I_m$$
удовлетворяет предельному соотношению
$$\|Z_n(\cdot)-I_m\|_\infty \rightarrow 0, \quad n\rightarrow \infty,$$
где $I_m$ --- единичная $(m\times m)$-матрица.

Положим теперь
$$A_{F_n}(\cdot):=\left(
                    \begin{matrix}
                      A_n(\cdot) & F_n(\cdot) \\
                      O_m & O_m \\
                    \end{matrix}
                  \right)\in L([a,b];\mathbb{C}^{2m\times 2m}),
$$
$$R_{A_nF_n}(\cdot):=A_{F_n}(\cdot)-A_F(\cdot)\in L([a,b];\mathbb{C}^{2m\times 2m}),$$
где $O_m $ --- нулевая $(m\times m)$-матрица.

Основным результатом данной работы является

\begin{theorem}
В формулировке теоремы Кигурадзе можно заменить условия $\rm (I)$, $\rm (II)$ на одно более общее условие
\begin{equation}\label{6}
R_{A_n}(\cdot)\in \mathcal{M}^{m},
\end{equation}
а условие $\rm (V)$ заменить на
\begin{equation}\label{7}
R_{A_nF_n}(\cdot)\in \mathcal{M}^{2m}.
\end{equation}
\end{theorem}

Условия \eqref{6}, \eqref{7} не являются конструктивными поскольку отсутствуют описания классов $\mathcal{M}^{m}$ и $\mathcal{M}^{2m}$. Однако из результатов работ \cite{Levin-A,Levin,Nguen,Kodliuk} вытекают удобные для применений достаточные условия принадлежности последовательности матриц-функций к этому классу. Поэтому из теоремы 1 вытекает ряд утверждений, которые обобщают или дополняют теорему Кигурадзе и выражаются в явном виде.

\begin{theorem}\label{theo2} В формулировке теоремы Кигурадзе можно заменить условие $\rm (II)$ на более общее условие
\begin{enumerate}
\item[$\rm (II^*)$] $\|R_{A_n}(\cdot)R^{\vee}_{A_n}(\cdot)\|_1\rightarrow 0$, $n\rightarrow\infty,$
\end{enumerate}
и добавить условие
\begin{enumerate}
\item[$\rm (VI^{\ast})$] $\|R_{A_n}(\cdot)R^{\vee}_{F_n}(\cdot)\|_1\rightarrow 0$, $n\rightarrow\infty.$
\end{enumerate}
\end{theorem}

 Преимущества теоремы~\ref{theo2} перед теоремой Кигурадзе становятся более заметными, если рассмотреть их приложения к системам линейных дифференциальных уравнений порядка $r\geq 2$ вида
\begin{equation}
\label{8}
 y^{(r)}(t)+A_{r-1}(t)y^{(r-1)}(t)+\dots+A_0(t)y(t)=f(t)
 \end{equation}
 с общими неоднородными краевыми условиями вида
 \begin{equation}
 \label{9}
 B_jy=c_j,  \quad j\in \{1,2,\ldots,r\}=:[r].
\end{equation}
Здесь линейные непрерывные операторы $B_j\colon C^{(r-1)}([a,b];\mathbb{C}^{m})\rightarrow \mathbb{C}^{m},$ а матрицы-функции $A_{j-1}(\cdot),$ вектор-функция $f(\cdot)$ и вектор $c_j$ такие же, как в \eqref{r1} и \eqref{r2}. Пусть теперь наряду с \eqref{8}--\eqref{9} задана последовательность краевых задач
\begin{equation}
\label{10}
 y^{(r)}_n(t)+A_{r-1,n}(t)y^{(r-1)}_n(t)+\dots+A_{0,n}(t)y_n(t)=f_n(t),
 \end{equation}
 \begin{equation}
 \label{11}
 B_{j,n}y_n=c_{j,n},  \quad j\in[r], \quad n\in \mathbb{N}.
\end{equation}

Каждую из этих задач очевидным образом можно свести к общей неоднородной краевой задаче для системы уравнений  первого порядка. Применительно к этим задачам теорема Кигурадзе приобретает следующий вид.

\begin{theorem}\label{theo3} Пусть выполнены условия $(0)$ и при $j\in[r]$, $n\rightarrow \infty$,
\begin{enumerate}
\item[$\rm (I^{\prime})$]  $\|R_{A_{j-1,n}}^{\vee}(\cdot)\|_\infty \rightarrow 0$,

\item[$\rm (II^{\prime})$]  $\|R_{A_{j-1,n}}(\cdot)\|_\infty=O(1)$,

\item[$\rm (III^{\prime})$]   $B_{j,n}y\rightarrow B_jy$, $y\in C^{(r-1)}([a,b];\mathbb{C}^m)$.
\end{enumerate}
Тогда для достаточно больших $n$ \emph{задача} \eqref{10}--\eqref{11} однозначно всюду разрешима.

Если кроме того
\begin{enumerate}
\item[$\rm (IV^{\prime})$] $c_{j,n}\rightarrow c_j$,  $j\in[r]$,  $n\rightarrow \infty$,
\item[$\rm (V)$] $\|R_{F_n}^{\vee}(\cdot)\|_\infty \rightarrow 0$,  $n\rightarrow \infty$,
\end{enumerate}
то единственные решения краевой задачи \eqref{10}--\eqref{11} удовлетворяют соотношениям
\begin{equation*}
\|y^{(j-1)}(\cdot)-y^{(j-1)}_n(\cdot)\|_\infty \rightarrow 0, \quad n\rightarrow \infty, \quad j\in[r].
\end{equation*}
\end{theorem}

Из теоремы~\ref{theo2} в этом случае вытекает, что справедлива

\begin{theorem}\label{theo4} В формулировке теоремы~{\rm \ref{theo3}} можно заменить условие $(\rm II^{\prime})$ на
\begin{enumerate}
\item[$\rm (II^{**})$] $\|R_{A_{r-1,n}}(\cdot)R^{\vee}_{A_{j-1,n}}(\cdot)\|_1\rightarrow 0$,  $n\rightarrow\infty$, $j\in[r]$.
\end{enumerate}
Если добавить условие
\begin{enumerate}
\item[$\rm (VI^{\ast \ast})$] $\|R_{A_{r-1,n}}(\cdot)R^{\vee}_{F_n}(\cdot)\|_1\rightarrow 0$,  $n\rightarrow\infty$,
\end{enumerate}
которое заведомо выполнено, коль выполнены условия $\rm (II^{\prime})$ и $\rm (V^{\prime})$.
\end{theorem}

Отметим, что условия (II$^{**}$) и (VI$^{**}$) заведомо выполнены если $\|R_{A_{r-1,n}}(\cdot)\|_1=O(1)$. При этом нет никаких ограничений на последовательность $\{\|A_{A_{j-1,n}}(\cdot)\|_1\colon n\geq 1 \}$ при $j\in[r-1]$.

\renewcommand\refname{References}


\begin{thebibliography}{99}

\bibitem{Reid}
Reid W. T. Some limit theorems for ordinary differential systems // J. Diff. Equat.~--- 1967.~--- v.~3, №~3.~--- P.~423--439. (DOI: 10.1016/0022-0396(67)90042-3).


\bibitem{Opial} Opial Z. Continuous parameter dependence in linear systems of differential equations // Ibid ---  1967.~---  v.~3.~---  Р.~571--579. (DOI: 10.1016/0022-0396(67)90017-4).

\bibitem{Levin-A} Levin  А. Yu. Passage to the limit for nonsingular systems $\dot{X}=A_n (t)X$ // Sov. Math. Dokl.~---  1967.~---  v.~\textbf{176},  №~4.~--  С.~774~-- 777.

\bibitem{Levin} Levin  А. Yu. Problems of the theory of ordinary differential equations. I // Vestn. Yaroslav. Univ.~---  1973.~---  Issue~5.~---  С.~105--132. (Russian).

\bibitem{Nguen} Nguyen Tkhe Hoan. Dependence of the solutions of a linear system of differential equations on a parameter //~ Differential Equations.~---  1993.~---  V.~\textbf{29}, №~6.~--  P.~830--835.

\bibitem{Kig} Kiguradze I. T. Some singular boundary value problems for ordinary differential equations.~---  Tbilisi: Izdat. Tbilis. Univ., 1975.~---  352~pp. (Russian).

\bibitem{Kigu} Kiguradze I. T. Boundary-value problems for systems of ordinary differential equations~// J.~Soviet Math.~--- 1988.~--- V.~\textbf{43}, №~2.~--- P.~2259--2339.

\bibitem{Kodliuk}
Kodliuk T. I., Mikhailets V. A., Reva N. V. Limit theorems for one-dimensional boundary-value problems // Ukrainian Math. J.~--- 2013.~--- V.~\textbf{65}.~--- №~1.~--- P.~77--90. (DOI: 10.1007/s11253-013-0766-x).

\bibitem{Gnyp}
Gnyp E. V., Kodliuk T. I., Mikhailets V. A. Fredholm boundary-value problems with parameter in Sobolev spaces // Ukrainian Math. J.~--- 2015~--- V.~\textbf{67}, №~5~--- P.~658--667. (DOI: 10.1007/s11253-015-1105-1).

\bibitem{Mikhailets}
Mikhailets V. A., Murach A. A., Soldatov V. A. Continuity in a parameter of solutions to generic boundary-value problems // Electron. J. Qual. Theory Differ. Equat.~--- 2016.~--- №~87.~--- P.~1--16.

\end{thebibliography}
\end{document}